\documentclass[review]{elsarticle}

\usepackage{lineno,hyperref}
\modulolinenumbers[5]

\usepackage{hyperref}
\hypersetup{nesting=true,debug=true,naturalnames=true}
\usepackage{graphicx,amsmath,amsfonts,amssymb,upref}
\usepackage{tikz}
\usepackage{multirow}
\usepackage{verbatim}
\usepackage{url}
\usepackage{float}
\floatstyle{plaintop}
\restylefloat{table}
\allowdisplaybreaks[1]

\let\<\langle
\let\>\rangle

\let\uml\"

\newtheorem{theorem}{Theorem}[section]
\newtheorem{lemma}[theorem]{Lemma}

\newtheorem{proposition}[theorem]{Proposition}
\newproof{pf}{Proof}

\journal{Journal of \LaTeX\ Templates}

\newcommand{\C}{\mathbb{C}}

\newcommand{\Z}{\mathbb{Z}}

\begin{document}

\begin{frontmatter}

\title{An explicit formula for the $A$-polynomial of the knot with Conway's notation $C(2n, 4)$}

\author{Ji-Young Ham\fnref{1}}
\address{School of Liberal Arts, Seoul National University of Science and Technology,
232 Gongneung-ro, Nowon-gu, Seoul, korea 
01811 \\
Department of Architecture, Konkuk University, 
120 Neungdong-ro, Gwangjin-gu, Seoul, Korea
05029
}
\ead{jiyoungham1@gmail.com} 
\fntext[1]{This work was supported by Basic Science Research Program through the National Research Foundation of Korea (NRF) funded by the Ministry of Education, Science and Technology (No. NRF-2018R1A2B6005847). }

\author{Joongul Lee\fnref{2}}  
\address{Department of Mathematics Education, Hongik University, 
94 Wausan-ro, Mapo-gu, Seoul, Korea
04066} 
\ead{jglee@hongik.ac.kr}  
\fntext[2]{The author was supported by 2018 Hongik University Research Fund.}

\begin{abstract}
An explicit formula for the $A$-polynomial of the knot having Conway's notation $C(2n,4)$ is computed  up to repeated factors. Our polynomial contains exactly the same irreducible factors as the $A$-polynomial defined in~\cite{CCGLS1}.
\end{abstract}

\begin{keyword}
$A$-polynomial, explicit formula, knot with Conway's notation $C(2n,4)$, Riley-Mednykh polynomial
\end{keyword}

\end{frontmatter}

\linenumbers

\section{Introduction}
In 1994, Cooper, Culler, Gillet, Long, and Shalen introduced the $A$-polynomial, $A(L,M) \in Z[L,M]$, of a compact 3-manifold $N$ with a single torus boundary~\cite{CCGLS1}. A-polynomials are related to the following invariants: incompressible surfaces~\cite{CCGLS1}, the Culler-Shalen seminorm~\cite{BZ1, CGLS}, cusp shapes~\cite{CL}, the volumes~\cite{CCGLS1,HL1,HMP} and Chern-Simons invariants~\cite{CS,HL,HL2,HL4,HLM3} of the deformed orbifolds, Mahler measure~\cite{Boyd}, Alexander polynomials~\cite{CCGLS1}, $AJ$ conjecture~\cite{FGL,Garou}, etc. The $AJ$ conjecture has been proved for some knots. For example, our knot, the knot with Conway's notation $C(2n,4)$ has been proved to satisfy the $AJ$ conjecture~\cite{LT} except $C(-4,4)$. 

 $A$-polynomials are known only for a few due to the difficulty of computations. Among them some are not full ones. A-polynomials were obtained using triangulation in~\cite{Champ,HMJP,HMJPT,Z2,Z3}. In~\cite{Culler}, Culler presented $A$-polynomials up to $8$ crossings and most $9$ crossings and many $10$ crossings, and all knots that can be triangulated with up to seven ideal tetrahedra. $A$-polynomials are known for twist knots~\cite{HS}, 
$(-2,3,1+2n)$ pretzel knots~\cite{GT,TY}, $J[m,2n]$~\cite[for $m$ between $2$ and $3$]{HS} ~\cite[for $m$ between $2$ and $5$ and for $m=2n$]{Pe}, recursively. Explicit formulas for the A-polynomials were computed for two-bridge torus knots~\cite{CCGLS1,HS,Fumikazu05}, iterated torus knots~\cite{Ni}, twist knots~\cite{HL,Mat1, Mat,Thompson}, knots with Conway’s notation $C(2n,3)$~\cite{HL3}, and the twisted torus knots $T (5, 1- 5n, 2, 2)$~\cite{Thompson}. We recall here that $J[4,-2n]$ is the mirror image of $C(2n,4)$.
\medskip

The main purpose of the paper is to find the explicit formula for the $A$-polynomial of the knot having Conway's notation $C(2n,4)$ up to repeated factors. Let us denote the $A$-polynomial of the knot having Conway's notation $C(2n,4)$ by $A_{2n}$. The following theorem gives the explicit formula for the $A$-polynomial of $C(2n,4)$.

\begin{theorem}  \label{thm:A-polynomial}
The $A$-polynomial $A_{2n}=A_{2n}(L,M)$ is given explicitly by 

\medskip
{\footnotesize \begin{align*}
 & \qquad \qquad \qquad A_{2n} =p_{2n}(u)p_{2n}(-u)\\
\text{where} \\
 p_{2n}(z)   = 
&  \begin{cases} 
 \sum_{i=0}^{2n} 
\binom{\left\lfloor \frac{i}{2}\right\rfloor +n}{i}
2^{-2 \left\lfloor \frac{i+1}{2}\right\rfloor -i} 
 \left(M^2\right)^{-\left\lfloor \frac{i}{2}\right\rfloor -2 \left\lfloor
   \frac{i+1}{2}\right\rfloor +i+n}
\left(L M^2+1\right)^{-2 \left\lfloor \frac{i+1}{2}\right\rfloor -i+2n}\\
\times \left(-2 L M^6+L M^4-L M^2-M^4+M^2 z+M^2-2\right)^{\left\lfloor
   \frac{i+1}{2}\right\rfloor } \\
\times \left(L M^2+L+M^2+z+1\right)^i \left(-3 L M^2+L+M^2+z-3\right)^{\left\lfloor
   \frac{i-1}{2}\right\rfloor } \\ 
\times \left((-1)^{i+1} \left(L M^2+1\right)-2 L M^2+L+M^2+z-2\right)        
\qquad \qquad 
\text{if $n \geq 0$,} \\
\sum_{i=0}^{-2n} 
\binom{\left\lfloor \frac{i-1}{2}\right\rfloor -n}{i}
2^{-2\left\lfloor \frac{i+1}{2}\right\rfloor -i} 
\left(M^2\right)^{-\left\lfloor \frac{i}{2}\right\rfloor -2 \left\lfloor
   \frac{i+1}{2}\right\rfloor +i-n}
\left(L M^2+1\right)^{-\frac{1}{2}-2 \left\lfloor \frac{i+1}{2}\right\rfloor -i-2n} \\
\times \left(-2 L M^6+L M^4-L M^2-M^4+M^2 z+M^2-2\right)^{\left\lfloor
   \frac{i+1}{2}\right\rfloor }\\
\times \left(L M^2+L+M^2+z+1\right)^i \left(-3 L M^2+L+M^2+z-3\right)^{\left\lfloor
   \frac{i-1}{2}\right\rfloor }\\       
\times \left((-1)^i \left(-2 L M^2+L+M^2+z-2\right)-L M^2-1\right)    
 \qquad \qquad 
\text{if $n<0$,}
\end{cases} \\
\text{and} \\
u &=\sqrt{5 L^2 M^4-2 L^2 M^2+L^2-2 L M^4+12 L M^2-2 L+M^4-2 M^2+5}.
\end{align*}}
\end{theorem}

\medskip 

Our writing is parallel with that in~\cite{HL3} which is based on~\cite{HL,HS,Mat}. The $A$-polynomial of $C(2n,4)$ can be obtained from the Riley-Mednykh polynomial in~\cite{HLMR}. With our theorem, we can easily get $A$-polynomials for $C(2n,4)$ having more than 12 crossings.
\section{Proof of Theorem~\ref{thm:A-polynomial}}

\begin{figure} 
\begin{center}
\resizebox{3.5cm}{!}{\includegraphics{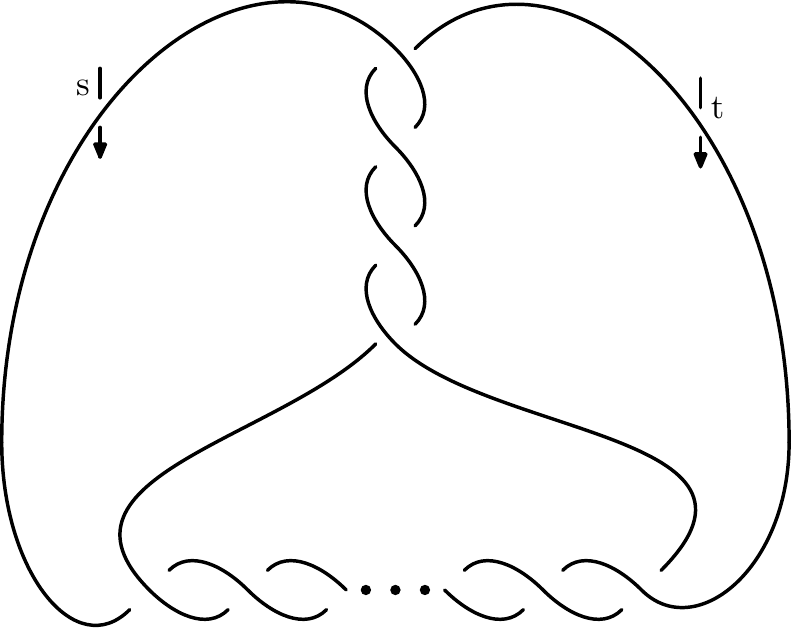}}
\qquad
\resizebox{3.5cm}{!}{\includegraphics{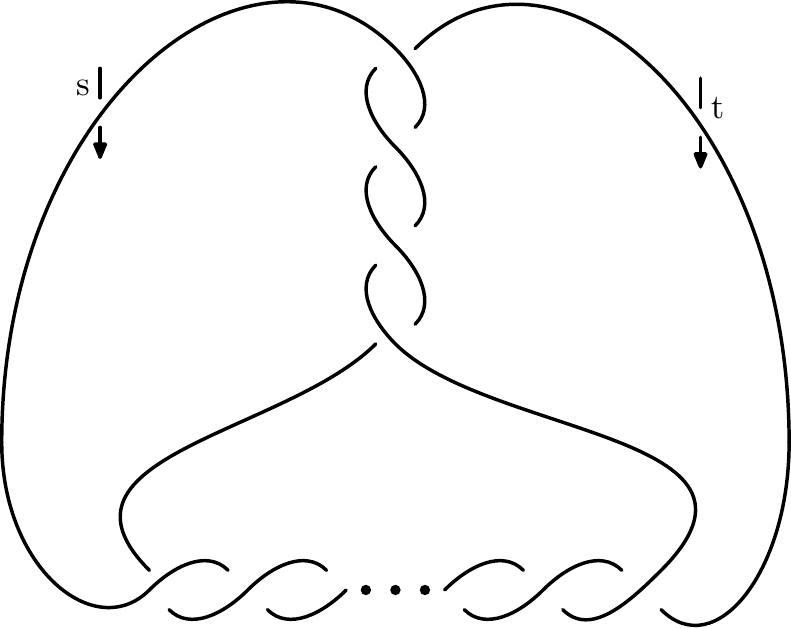}}
\caption{A two bridge knot having Conway's notation $C(2n,4)$  for $n>0$ (left) and for $n<0$ (right)} \label{fig:C[2n,4]}
\end{center}
\end{figure}

A knot $K$ is a two bridge knot with Conway's notation $C(2n,4)$ if $K$ has a regular two-dimensional projection of the form in Figure~\ref{fig:C[2n,4]}. Let us denote the exterior of $C(2n,4)$ by $X_{2 n}$. The following proposition gives the fundamental group of 
$X_{2n}$~\cite{HS,HLMR,R1}.
\begin{proposition}\label{theorem:fundamentalGroup}
$$\pi_1(X_{2n})=\left\langle s,t \ |\ swt^{-1}w^{-1}=1\right\rangle,$$
where $w=(ts^{-1}ts^{-1}t^{-1}st^{-1}s)^n$.
\end{proposition}
\medskip 

Given  a set of generators, $\{ s,t \}$, of the fundamental group for 
$\pi_1 (X_{2n})$, we define a representation $\rho \ : \ \pi_1 (X_{2n}) \rightarrow \text{SL}(2, \C)$ by

{\footnotesize \begin{center}
$$\begin{array}{ccccc}
\rho(s)=\left[\begin{array}{cc}
                       M &       1 \\
                        0      & M^{-1}  
                     \end{array} \right]                          
\text{,} \ \ \
\rho(t)=\left[\begin{array}{cc}
                   M &  0      \\
                  2-M^2-M^{-2} -x  \ \ \   & M^{-1} 
                 \end{array}  \right].
\end{array}$$
\end{center}}

Then $\rho$ can be identified with the point $(M,x) \in \C^2$.  When $M$ varies we have an algebraic set 
whose defining equation is the following explicit Riley-Mednykh polynomial $P_{2n}=P_{2n}(M,x)$ defined recursively by the following form: 

\begin{align} \label{equ:P}
P_{2n} = \begin{cases}
 Q P_{2(n-1)} -M^{12} P_{2(n-2)} \ \text{if $n>1$}, \\
 Q P_{2(n+1)}-M^{12} P_{2(n+2)} \ \text{if $n<-1$},
\end{cases}
\end{align}
\medskip

\noindent with initial conditions
{\footnotesize \begin{align*} 
P_{-2}  & =-M^4 x^3-\left(2 M^6-M^4+2 M^2\right) x^2-\left(M^8-M^6+2
   M^4-M^2+1\right) x+M^4,\\
P_{0}  & =M^{-2} \ \text{\textnormal{for}} \ n < 0 \qquad \text{\textnormal{and}} \qquad P_{0} =1 \ \text{\textnormal{for}} \ n \geq 0,\\    
P_{2}  & =M^6 x^4+\left(3 M^8-M^6+3 M^4\right) x^3+\left(3 M^{10}-2 M^8+5 M^6-2 M^4+3
   M^2\right) x^2\\
   &+\left(M^{12}-M^{10}+2 M^8-2 M^6+2 M^4-M^2+1\right) x+M^6, \\
\end{align*}
\noindent  and
\begin{align*}
Q&=M^6 x^4+\left(3 M^8-2 M^6+3 M^4\right) x^3+\left(3 M^{10}-4 M^8+6 M^6-4 M^4+3
   M^2\right) x^2\\
   &+\left(M^{12}-2 M^{10}+3 M^8-4 M^6+3 M^4-2 M^2+1\right) x+2 M^6.
\end{align*}}

It is known in~\cite{HLMR} that $\rho$ is a representation of $\pi_1(X_{2n})$ if and only if $x$ is a root of $P_{2n}$. We use different initial conditions $P_0$ for $P_{2n}(n>1)$ and $P_{2n}(n<-1)$, but the same $P$ for both series for simplicity.

\begin{lemma} \label{lem1:RMpolynomial}
The Riley-Mednykh polynomial $P_{2n}=P_{2n}(x,M)$ is described explicitly by 

\begin{align*} \footnotesize
P_{2n} = \begin{cases}
 \sum_{i=0}^{2n} 
\binom{\left\lfloor \frac{i}{2}\right\rfloor +n}{i}  \left(M^2\right)^{-\left\lfloor \frac{i}{2}\right\rfloor -2 \left\lfloor
\frac{i+1}{2}\right\rfloor +3 n} x^{\left\lfloor \frac{i+1}{2}\right\rfloor }
\left(M^4+M^2 x-2 M^2+1\right)^{\left\lfloor \frac{i-1}{2}\right\rfloor }   \\
\times  \left(M^4+M^2 x+1\right)^i 
\left(\left((-1)^{i+1}+1\right) M^2/2+M^4+M^2 x-2 M^2+1\right)  
\qquad \ \ \\ \text{if $n \geq 0$}, \\
 \sum_{i=0}^{-2n} 
\binom{\left\lfloor \frac{i-1}{2}\right\rfloor -n}{i} 
\left(M^2\right)^{-\left\lfloor \frac{i}{2}\right\rfloor -2 \left\lfloor
   \frac{i+1}{2}\right\rfloor -3 n-1}
x^{\left\lfloor \frac{i+1}{2}\right\rfloor }   
\left(M^4+M^2 x-2 M^2+1\right)^{\left\lfloor \frac{i-1}{2}\right\rfloor }  \\
\times \left(M^4+M^2 x+1\right)^i 
\left(\left((-1)^i-1\right) M^2/2+(-1)^i
\left(M^4+M^2 x-2 M^2+1\right)\right) 
\ \\ \text{if $n<0$}.
\end{cases}
\end{align*}
\end{lemma}

\begin{pf}
Write $f_{2n}$ for the claimed formula and show that $f_{2n}=P_{2n}$. 

\noindent Case I: $n \geq 0$. When $i > 2n$ or $i<0$, 
$\binom{\left\lfloor \frac{i}{2}\right\rfloor +n}{i}$
is undefined and can be considered as zero. Hence the finite sum can be regarded as an infinite sum.
Direct computation shows that $f_0=P_0$ and $f_2=P_2$.
Now, we only need to show that $f_{2n}$ satisfies the recursive relation.
Note that $Q$ can be written as $x \left(M^4+M^2 x+1\right)^2 \left(M^4+M^2 x-2 M^2+1\right)+2 M^6$.
{\footnotesize \begin{align*}
&Q f_{2(n-1)} -M^{12} f_{2(n-2)} \\
&=
 \left(x \left(M^4+M^2 x+1\right)^2 \left(M^4+M^2 x-2 M^2+1\right)+2 M^6\right) \\
& \times 
\sum_i
 \binom{\left\lfloor \frac{i}{2}\right\rfloor +n-1}{i}  \left(M^2\right)^{-\left\lfloor \frac{i}{2}\right\rfloor -2 \left\lfloor
\frac{i+1}{2}\right\rfloor +3 (n-1)} x^{\left\lfloor \frac{i+1}{2}\right\rfloor }
\left(M^4+M^2 x-2 M^2+1\right)^{\left\lfloor \frac{i-1}{2}\right\rfloor }
  \\
& \times  \left(M^4+M^2 x+1\right)^i 
\left(\left((-1)^{i+1}+1\right) M^2/2+M^4+M^2 x-2 M^2+1\right)  \\
&-M^{12}
\sum_i
\binom{\left\lfloor \frac{i}{2}\right\rfloor +n-2}{i}  \left(M^2\right)^{-\left\lfloor \frac{i}{2}\right\rfloor -2 \left\lfloor
\frac{i+1}{2}\right\rfloor +3 (n-2)} x^{\left\lfloor \frac{i+1}{2}\right\rfloor }
\left(M^4+M^2 x-2 M^2+1\right)^{\left\lfloor \frac{i-1}{2}\right\rfloor } \\
& \times  \left(M^4+M^2 x+1\right)^i 
\left(\left((-1)^{i+1}+1\right) M^2/2+M^4+M^2 x-2 M^2+1\right)  \\
& = \sum_i 
\left[
\binom{\left\lfloor \frac{i}{2}\right\rfloor +n-2}{i-2}
+2 \binom{\left\lfloor \frac{i}{2}\right\rfloor +n-1}{i}
-\binom{\left\lfloor \frac{i}{2}\right\rfloor +n-2}{i}
\right] \\
&   \left(M^2\right)^{-\left\lfloor \frac{i}{2}\right\rfloor -2 \left\lfloor
\frac{i+1}{2}\right\rfloor +3 n} x^{\left\lfloor \frac{i+1}{2}\right\rfloor }
\left(M^4+M^2 x-2 M^2+1\right)^{\left\lfloor \frac{i-1}{2}\right\rfloor } \\
& \times  \left(M^4+M^2 x+1\right)^i 
\left(\left((-1)^{i+1}+1\right) M^2/2+M^4+M^2 x-2 M^2+1\right) \\
&= \sum_i
\binom{\left\lfloor \frac{i}{2}\right\rfloor +n}{i}
\left(M^2\right)^{-\left\lfloor \frac{i}{2}\right\rfloor -2 \left\lfloor
\frac{i+1}{2}\right\rfloor +3 n} x^{\left\lfloor \frac{i+1}{2}\right\rfloor }
\left(M^4+M^2 x-2 M^2+1\right)^{\left\lfloor \frac{i-1}{2}\right\rfloor } \\
& \times  \left(M^4+M^2 x+1\right)^i 
\left(\left((-1)^{i+1}+1\right) M^2/2+M^4+M^2 x-2 M^2+1\right)\\
&=f_{2n}
\end{align*}}
In the last equality we use the binomial relation 
{\footnotesize $$\binom{a}{b} =\binom{a-1}{b-1}+\binom{a-1}{b}$$}
three times.

\noindent Case II: $n < 0$. When $i > -2n-1$ or $i<0$, 
$\binom{\left\lfloor \frac{i-1}{2}\right\rfloor -n}{i}$
is undefined and can be considered as zero. Hence the finite sum can be regarded as an infinite sum.
Direct computation shows that $f_0=P_0$ and $f_{-2}=P_{-2}$. As in Case I, one can easily show that $f_{2n}$ satisfies the recursive relation.
\qed
\end{pf}

   Let $l = ww^{*}$~\cite{CCGLS1,HS}, where $w^{*}$ is the word obtained by reversing $w$ (by reading words in $w$ from right to left). Let $L=\rho(l)_{11}$ (the left upper entry of $\rho(l)$, \cite[lemma 7.9]{HLMR}). Then $l$ is the longitude which is null-homologous in $X_{2n}$ (you can read a twisted longitude $ww^{*}$ from the Schubert normal form of the knot $C(2n,4)$), and we have
   
\medskip 
   
\begin{lemma}\cite[Theorem 7.10]{HLMR}
\label{lem2:longitude}
{\footnotesize \begin{align*}
L &=-M^{-2}\frac{M^{-4}-M^{-2}+(2 M^{-2}+M^{2}-1)x+x^2}{M^{4}-M^{2}+(M^{-2}+2M^{2}-1)x+x^2}, \\
x &=\frac{-2 L M^8+L M^6-L M^4-M^6 +M^4 z+M^4-2 M^2}{2 \left(L M^6+M^4\right)} \\
\text{where} \\
z  = u &= \sqrt{5 L^2 M^4-2 L^2 M^2+L^2-2 L M^4+12 L M^2-2 L+M^4-2 M^2+5} \\
\text{ or } z= -u.
\end{align*} }
\end{lemma}
Now substituting $\frac{-2 L M^8+L M^6-L M^4-M^6 +M^4 z+M^4-2 M^2}{2 \left(L M^6+M^4\right)}$ for $x$ into $P_{2n}$, for $n > 0$, gives 
{\tiny
\begin{align*} 
& \sum_{i=0}^{2n} 
\binom{\left\lfloor \frac{i}{2}\right\rfloor +n}{i}  \left(M^2\right)^{-\left\lfloor \frac{i}{2}\right\rfloor -2 \left\lfloor\frac{i+1}{2}\right\rfloor +3 n} 
 \left(\frac{-2 L M^8+L M^6-L M^4-M^6+ M^4 z+M^4-2 M^2}{2 \left(L M^6+M^4\right)} \right)^{\left\lfloor \frac{i+1}{2}\right\rfloor } \\
& \times\left(M^4+M^2 \frac{-2 L M^8+L M^6-L M^4-M^6 +M^4 z+M^4-2 M^2}{2 \left(L M^6+M^4\right)}-2 M^2+1\right)^{\left\lfloor \frac{i-1}{2}\right\rfloor }  \\
& \times  \left(M^4+M^2 \frac{-2 L M^8+L M^6-L M^4-M^6 +M^4 z+M^4-2 M^2}{2 \left(L M^6+M^4\right)}+1\right)^i \\
& \times \left(\left((-1)^{i+1}+1\right) M^2/2+M^4+M^2 \frac{-2 L M^8+L M^6-L M^4-M^6 +M^4 z+M^4-2 M^2}{2 \left(L M^6+M^4\right)}-2 M^2+1\right) 
\end{align*}
}

We observe that 
{\tiny
\begin{align*}
& \frac{-2 L M^8+L M^6-L M^4-M^6 +M^4 z+M^4-2 M^2}{2 \left(L M^6+M^4\right)}=\frac{M^2 \left(-2 L M^6+L M^4-L M^2-M^4 +M^2 z+M^2-2\right)}{2 M^4 \left(L M^2+1\right)},\\
& M^4+M^2 \frac{-2 L M^8+L M^6-L M^4-M^6 +M^4 z+M^4-2 M^2}{2 \left(L M^6+M^4\right)}-2 M^2+1= \frac{M^6 \left(-3 L M^2+L+M^2+z-3\right)}{2 M^4 \left(L M^2+1\right)}, \\
& M^4+M^2 \frac{-2 L M^8+L M^6-L M^4-M^6 +M^4 z+M^4-2 M^2}{2 \left(L M^6+M^4\right)}+1= \frac{M^6 \left(L M^2+L+M^2+z+1\right)}{2 M^4\left(L M^2+1\right)},\\
\end{align*}
}
and
{\tiny
\begin{align*}
& \left((-1)^{i+1}+1\right) M^2/2+M^4+M^2 \frac{-2 L M^8+L M^6-L M^4-M^6 +M^4 z+M^4-2 M^2}{2 \left(L M^6+M^4\right)}-2 M^2+1 \\
& =\frac{M^6 \left((-1)^{i+1} \left(L M^2+1\right)-2 L M^2+L+M^2+z-2\right)}{2 M^4 \left(L M^2+1\right)}.
\end{align*}
}

The resulting expression is

{\tiny
\begin{align*} 
q_{2n}(z)&= \sum_{i=0}^{2n} 
\binom{\left\lfloor \frac{i}{2}\right\rfloor +n}{i}  \left(M^2\right)^{-\left\lfloor \frac{i}{2}\right\rfloor -2 \left\lfloor\frac{i+1}{2}\right\rfloor +3 n} 
 \left(\frac{-2 L M^6+L M^4-L M^2-M^4 +M^2 z+M^2-2}{2 M^2 \left(L M^2+1\right)} \right)^{\left\lfloor \frac{i+1}{2}\right\rfloor } \\
& \times\left(\frac{M^2 \left(-3 L M^2+L+M^2+z-3\right)}{2  \left(L M^2+1\right)}\right)^{\left\lfloor \frac{i-1}{2}\right\rfloor }  
 \left(\frac{M^2 \left(L M^2+L+M^2+z+1\right)}{2 \left(L M^2+1\right)}\right)^i \\
& \times \left(\frac{M^2 \left((-1)^{i+1} \left(L M^2+1\right)-2 L M^2+L+M^2+z-2\right)}{2 \left(L M^2+1\right)}\right). 
\end{align*}
}

Since $q_{2n}(u)$ is a non-polynomial factor of the $A$-polynomial, we multiply it by its Galois conjugate $q_{2n}(-u)$ to obtain the entire $A$-polynomial~\cite{CCGLS1,HS}. This is actually another method of defining the resultant of two polynomials~\cite[p.209-p.210]{Lang}.
We multiply $q_{2n}(u) q_{2n}(-u)$ by 
$M^{-8n} \left(L M^2+1\right)^{4n}$ to clear the denominators and factor out some power of $M$ so that we have the $A$-polynomial $ A_{2n}(L,M)=p_{2n}(u)p_{2n}(-u)$  in Theorem~\ref{thm:A-polynomial}. Formally, it can be done by multiplying $q_{2n}(z)$ by $M^{-4n} \left(L M^2+1\right)^{2n}$. Notice that $p_{2n}(z)=M^{-4n} \left(L M^2+1\right)^{2n} q_{2n}(z)$ and recall that 
{\tiny $$u =\sqrt{5 L^2 M^4-2 L^2 M^2+L^2-2 L M^4+12 L M^2-2 L+M^4-2 M^2+5}.$$ }
\noindent Now we want to show that the claimed formula does not have fractions.
Direct computation shows that  $A_0(L,M)=1$, $A_2(L,M)$, and $A_4(L,M)$ (see, Appendix B) are polynomials in $\Z[L,M]$. From now on a polynomial means a polynomial in $\Z[L,M]$. Let us assume that $A_{2k}$ for $k \leq n$ is a polynomial.
From the equation~(\ref{equ:P}), we have $P_{2n}=Q P_{2(n-1)} -M^{12} P_{2(n-2)}$ for $n>1$. Hence 
$q_{2n}=Q q_{2(n-1)} -M^{12} q_{2(n-2)}$
and 
{\footnotesize \begin{align} \label{equ1}
p_{2n}(z)=M^{-4} \left(L M^2+1\right)^{2} Q p_{2(n-1)}(z)
- \left(L M^2+1\right)^{4} M^{4} p_{2(n-2)} (z). 
\end{align}}
Now, we have 
\begin{proposition} \label{pro1}
{\footnotesize \begin{align}
&A_{2n}=p_{2n}(u) p_{2n}(-u) \nonumber \\
&=M^{-8} \left(L M^2+1\right)^{4} Q(u) Q(-u) p_{2(n-1)}(u) p_{2(n-1)}(-u)  \label{pro1-1}\\
&- \left(L M^2+1\right)^{6} \left(Q(u) p_{2(n-1)}(u)p_{2(n-2)}(-u)+Q(-u) p_{2(n-1)}(-u)p_{2(n-2)}(u) \right) \label{pro1-2}\\
&+M^{8} \left(L M^2+1\right)^{8} p_{2(n-2)} (u) p_{2(n-2)} (-u). \label{pro1-3}
\end{align}}
\end{proposition}
By direct computations, 
{\footnotesize \begin{align}
Q(u) &=2^{-1} M^4 \left(L M^2+1\right)^{-4} \left(a(L,M)+L b(L,M) u \right) \label{equq-1}\\ 
 Q(-u) &=2^{-1} M^4 \left(L M^2+1\right)^{-4} \left(a(L,M)-L b(L,M) u\right) \label{equq-2}\\
 Q(u) Q(-u) &= M^{8} \left(L M^2+1\right)^{-4} f(L,M) \label{equqq} \\
 Q(u) Q(u) &=2^{-1} M^{8} \left(L M^2+1\right)^{-8}\left( g(L,M)+h(L,M) u\right) \label{equqqp} \\
Q(-u) Q(-u) &=2^{-1} M^{8} \left(L M^2+1\right)^{-8} \left(g(L,M)-h(L,M) u\right) \label{equqqm}
\end{align} }
 where $a(L,M)$, $b(L,M)$, $f(L,M)$, $g(L,M)$ and $h(L,M)$ are some polynomials (see Appendix A).
 Since $A_{2n}=p_{2n}(u) p_{2n}(-u)$, $A_{2(n-1)}=p_{2(n-1)}(u) p_{2(n-1)}(-u)$, and $$A_{2(n-2)}=p_{2(n-2)}(u) p_{2(n-2)}(-u)$$ are polynomials and since $M^{-8} \left(L M^2+1\right)^{4} Q(u) Q(-u)$ is a polynomial  by Equation~\eqref{equqq}, the terms in 
 Proposition~\ref{pro1}-\eqref{pro1-1} and Proposition~\ref{pro1}-\eqref{pro1-3} are polynomials and hence the term
{\footnotesize $$-\left(L M^2+1\right)^{6} \left(Q(u) p_{2(n-1)}(u)p_{2(n-2)}(-u)+Q(-u) p_{2(n-1)}(-u)p_{2(n-2)}(u) \right)$$} in Proposition~\ref{pro1}-\eqref{pro1-2} is a polynomial.
 Now we only need to show that $A_{2(n+1)}=p_{2(n+1)}(u) p_{2(n+1)}(-u)$ is a polynomial. By Proposition~\ref{pro1}, we have
{\footnotesize \begin{align}
& A_{2(n+1)}=p_{2(n+1)}(u) p_{2(n+1)}(-u)  \nonumber \\
&=M^{-8} \left(L M^2+1\right)^{4} Q(u) Q(-u) p_{2n}(u) p_{2n}(-u)  \label{pro1-4}\\
&- \left(L M^2+1\right)^{6} \left(Q(u) p_{2n}(u)p_{2(n-1)}(-u)+Q(-u) p_{2n}(-u)p_{2(n-1)}(u) \right) \label{pro1-5}\\
&+M^{8} \left(L M^2+1\right)^{8} p_{2(n-1)} (u) p_{2(n-1)} (-u). \label{pro1-6}
 \end{align}}
 Since the terms in~\eqref{pro1-4} and~\eqref{pro1-6} are polynomials by the induction hypothesis and  
 Equation~\eqref{equqq}, we only need to show the term in~\eqref{pro1-5} is a polynomial.
 Now $-1$ times the term in \eqref{pro1-5} becomes 
{\footnotesize \begin{align}
 & M^{-4}\left(L M^2+1\right)^{8} \left(Q(u)  Q(u)+Q(-u)  Q(-u)\right) p_{2(n-1)}(u)p_{2(n-1)}(-u)  \label{equ2-1}\\
&- M^{4}\left(L M^2+1\right)^{10}  (Q(u) p_{2(n-1)}(-u) p_{2(n-2)} (u)+Q(-u) p_{2(n-1)}(u) p_{2(n-2)} (-u))  \label{equ2-2}
  \end{align} }
 by using Equation~\eqref{equ1}.
 The term in~\eqref{equ2-1} is a polynomial because of Equations~\eqref{equqqp},~\eqref{equqqm}, and the induction hypothesis.
 We will show that the term in~\eqref{equ2-2} is a polynomial.
 By using Equation~\eqref{equ1} again, $-1$ times the term in~\eqref{equ2-2} becomes

{\footnotesize \begin{align*}
&2 \left(L M^2+1\right)^{12} Q(u) Q(-u) p_{2(n-2)}(-u) p_{2(n-2)} (u) \\
&- M^{8} \left(L M^2+1\right)^{14} \left(Q(u) p_{2(n-2)} (u) p_{2(n-3)} (-u) +Q(-u) p_{2(n-2)} (-u) p_{2(n-3)} (u) \right).   \\
 \end{align*}}
 It is a polynomial by~\eqref{pro1-2},~\eqref{equqq}, and the induction hypothesis.
Therefore, for all $n>0$, by induction, $ A_{2n}(L,M)=p_{2n}(u)p_{2n}(-u)$ is a polynomial. 

Now, we are going to show that $A_{2n}(L,M)$ for $n>0$ has $L^2$ and $M^{8n}$ as terms so that it does not have any redundant $L$ or $M$ factors.
Let us consider $p_{2n}$ and $Q$ as functions of $L$, $M$ and $z$ and $u$ as a function of $L$ and $M$. Let 
$(\beta_1)=(L,0,u(L,0))$ and $(\bar{\beta_1})=(L,0,-u(L,0))$.
When $M$ is equal to zero, Proposition~\ref{pro1} simplifies to
\begin{proposition} \label{pro2}
 \begin{align*} 
&p_{2n}(\beta_1) p_{2n}(\bar{\beta_1})\\
&=f(L,0) p_{2(n-1)}(\beta_1) p_{2(n-1)}(\bar{\beta_1}) \\
&=(1-3L+3L^2-L^3) p_{2(n-1)}(\beta_1) p_{2(n-1)}(\bar{\beta_1}).
\end{align*} 
\end{proposition}
Since 
\begin{align}
p_{2}(\beta_1) p_{2}(\bar{\beta_1})=L^2-L^3, \label{equb}
\end{align}
 for all $n>0$, by induction, $ A_{2n}(L,M)=p_{2n}(u)p_{2n}(-u)$ has $L^2$ as a term.
Let $(\beta_2)=\left(0,M,u(0,M)\right)$ and $(\bar{\beta_2})=\left(0,M,-u(0,M)\right)$. 
Since $Q(\beta_2)=Q(\bar{\beta_2})=M^4 (1+M^4)$, we have the following equations from Equation~(\ref{equ1}):
\begin{align*} 
p_{2n}(\beta_2)&=(1+M^4) p_{2(n-1)}(\beta_2)
-  M^{4} p_{2(n-2)} (\beta_2) \\
p_{2n}(\bar{\beta_2})&=(1+M^4) p_{2(n-1)}(\bar{\beta_2})
-  M^{4} p_{2(n-2)} (\bar{\beta_2}). 
\end{align*}
Since $p_0=1$, $p_2(\beta_2)=p_2(\bar{\beta_2})=M^4$, by induction, we have $p_{2n}(\beta_2)=p_{2n}(\bar{\beta_2})=M^{4n}$ and $A_{2n}(0,M)=M^{8n}$. Therefore $A_{2n}(L,M)$ for $n>0$ has $M^{8n}$ as a term.
Now, we are going to show that $A_{2n}(L,M)$ for $n>0$ has $L^{3n}$ as a term and does not have a constant (nonzero) times $L^{3n+1} M^{2}$  as a term so that
it does not have any redundant $LM^2+1$ factors.
By Proposition~\ref{pro2} and Equation~\eqref{equb}, $A_{2n}(L,M)$ for $n>0$ has $L^{3n}$ as a term.
Since $p_0=1$ and $p_2(\pm u)$ (see Appendix C) is a polynomial up to $\left(2 \left(L M^2+1\right)\right)^{2}$, by Equations~\eqref{equ1}, ~\eqref{equq-1} and~\eqref{equq-2}, $p_{2 n}(\pm u)$ for $n>0$ is a polynomial up to $\left(2 \left(L M^2+1\right)\right)^{2 n}$. Hence, all the terms in~(\ref{pro1-3}) have at least 8th power of $M$, and, by using Equations~\eqref{equq-1},~\eqref{equq-2} again, all the terms in~(\ref{pro1-2}) have at least 4th power of $M$ because we know that factoring out $\left(2 \left(L M^2+1\right)\right)^{*}$ from a polynomial won't affect $M^{**}$ factor of it and formulae in~(\ref{pro1-2}) and (\ref{pro1-3}) are polynomials. Therefore if there were a constant (nonzero) times $L^{3n+1} M^{2}$ in $A_{2n}(L,M)$ for $n>0$, it would be in the terms in~(\ref{pro1-1}). 
By using Equation~(\ref{equqq}), the terms in~(\ref{pro1-1}) can be rewritten as 
\begin{align}
f(L,M) p_{2(n-1)}(u) p_{2(n-1)}(-u)  \label{pro1-1-1}.
\end{align}
The terms in $f(L,M)$ whose power of $M$ is less than or equal to $2$ are 
$$-L^3+3 L^2-3 L+1+\left(2 L^3-8 L^2+6 L\right) M^2.$$
The terms in $A_2(L,M)=p_2(u) p_2(-u)$ whose power of $M$ is less than or equal to $2$ are
$$L^3+L^2+\left(L^3+L^2\right) M^2.$$
Now, one can easily prove that the terms in~\ref{pro1-1-1} do not include a constant (nonzero) times $L^{3n+1} M^{2}$.

From now on, we will deal with the case for $n<0$.

 Substituting $\frac{-2 L M^8+L M^6-L M^4-M^6 +M^4 z+M^4-2 M^2}{2 \left(L M^6+M^4\right)}$ for $x$ into $P_{2n}$, for $n < 0$, gives
{\tiny
\begin{align*} 
q_{2n}(z)&= \sum_{i=0}^{-2n} 
\binom{\left\lfloor \frac{i-1}{2}\right\rfloor -n}{i} 
\left(M^2\right)^{-\left\lfloor \frac{i}{2}\right\rfloor -2 \left\lfloor
   \frac{i+1}{2}\right\rfloor -3 n-1}
 \left(\frac{-2 L M^6+L M^4-L M^2-M^4 +M^2 z+M^2-2}{2 M^2 \left(L M^2+1\right)} \right)^{\left\lfloor \frac{i+1}{2}\right\rfloor } \\
& \times\left(\frac{M^2 \left(-3 L M^2+L+M^2+z-3\right)}{2  \left(L M^2+1\right)}\right)^{\left\lfloor \frac{i-1}{2}\right\rfloor }  
 \left(\frac{M^2 \left(L M^2+L+M^2+z+1\right)}{2 \left(L M^2+1\right)}\right)^i \\
& \times \left(\frac{M^2 \left((-1)^i \left(-2 L M^2+L+M^2+z-2\right)-L M^2-1\right)}{2 \left(L M^2+1\right)}\right). 
\end{align*}
}

As in case $n>0$, we multiply it by its Galois conjugate $q_{2n}(-u)$ to obtain the entire $A$-polynomial.
We multiply $q_{2n}(u) q_{2n}(-u)$ by 
$M^{8n+4} \left(L M^2+1\right)^{-4n+1}$ to clear the denominators and factor out some power of $M$ so that we have the $A$-polynomial $ A_{2n}(L,M)=p_{2n}(u)p_{2n}(-u)$  in Theorem~\ref{thm:A-polynomial}. Similarly as in case $n>0$, we can show that $A_{2n}(L,M)$ for $n<0$ does not have any redundant $L$, $M$ or $LM^2+1$ factors. It can be done by showing that it has $1$ and $L^{3(n-1)+1}$ as terms and does not have a constant times $L^{3(n-1)+2} M^{2}$  as a term.

\section{Appendix A} \label{app:A}
 {\tiny
 \begin{align*}
& a(L,M) =2 L^4\left(M^{12} + M^8 \right) + 
 L^3\left(-3 M^{12} + 10 M^{10} + 5 M^8 + 4 M^6 - M^4 + 2 M^2 - 
    1 \right) +2 \left(M^4+1\right)\\
&+2 L^2 \left(M^{12}-4 M^{10}+5 M^8+8 M^6+5 M^4-4 M^2+1\right)-L
   \left(M^{12}-2 M^{10}+M^8-4 M^6-5 M^4-10 M^2+3\right) \\
& b(L,M) = (L-1) (M-1)^3 (M+1)^3 \left(M^2+1\right)^2\\
& f(L,M) = -L^3+3 L^2-3 L+1+\left(6 L^3-8 L^2+2 L\right) M^{14}+\left(6 L^3+8 L^2+2 L\right) M^{10}+\left(2 L^3+8 L^2+6 L\right) M^6 \\
  & +\left(-2 L^3+2 L^2+2 L+2\right) M^4+\left(2 L^3-8 L^2+6 L\right) M^2+\left(L^4-3 L^3+3 L^2-L\right) M^{16} \\
  & +\left(2 L^4+2 L^3+2 L^2-2 L\right) M^{12}+\left(L^4+4 L^3+14 L^2+4
   L+1\right) M^8. \\ 
& g(L,M) =2 L^8 \left(M^4+1\right)^2 M^{16}+2 L^7 \left(-3 M^{16}+10 M^{14}+2 M^{12}+14
   M^{10}+4 M^8+6 M^6-2 M^4+2 M^2-1\right) M^8 \\
&+L^6 \left(11 M^{24}-52 M^{22}+54
   M^{20}+76 M^{18}+91 M^{16}-8 M^{14}+68 M^{12}-8 M^{10}-7 M^8-4 M^6+6 M^4-4
   M^2+1\right) \\
& -2 L^5 \left(7 M^{24}-34 M^{22}+52 M^{20}+18 M^{18}-106 M^{16}-178
   M^{14}+54 M^{12}-22 M^{10}-31 M^8+4 M^6+22 M^4-12 M^2+2\right)\\
& +2 L^4 \left(5
   M^{24}-28 M^{22}+44 M^{20}+4 M^{18}-71 M^{16}+24 M^{14}+324 M^{12}+24 M^{10}-71
   M^8+4 M^6+44 M^4-28 M^2+5\right)\\
& -2 L^3 \left(2 M^{24}-12 M^{22}+22 M^{20}+4
   M^{18}-31 M^{16}-22 M^{14}+54 M^{12}-178 M^{10}-106 M^8+18 M^6+52 M^4-34
   M^2+7\right)\\
&+L^2 \left(M^{24}-4 M^{22}+6 M^{20}-4 M^{18}-7 M^{16}-8 M^{14}+68
   M^{12}-8 M^{10}+91 M^8+76 M^6+54 M^4-52 M^2+11\right) \\
& +L \left(-2 M^{16}+4
   M^{14}-4 M^{12}+12 M^{10}+8 M^8+28 M^6+4 M^4+20 M^2-6\right)+2
   \left(M^4+1\right)^2\\
& h(L,M) = L \left(L-1\right) \left(M^2-1\right)^3 \left(M^2+1\right)^2 h_1(L,M) \\
& h_1(L,M) =L^4\left (2 M^{12} + 2 M^8 \right) + 
 L^3\left(-3 M^{12} + 10 M^{10} + 5 M^8 + 4 M^6 - M^4 + 2 M^2 - 
    1 \right) \\
   & +L^2 \left(2 M^{12}-8 M^{10}+10 M^8+16 M^6+10 M^4-8 M^2+2\right)\\
   & +L
   \left(-M^{12}+2 M^{10}-M^8+4 M^6+5 M^4+10 M^2-3\right)+2 M^4+2
 \end{align*}
 }
\section{Appendix B} \label{app:B}
  {\tiny
 \begin{align*}
& A_2(L,M) =L^4 M^8+L^3 \left(-2 M^{12}+3 M^{10}+3 M^8+M^2-1\right)+L^2 \left(M^{16}-3
   M^{14}-M^{12}+3 M^{10}+6 M^8+3 M^6-M^4-3 M^2+1\right) \\
&   +L \left(-M^{16}+M^{14}+3 M^8+3 M^6-2 M^4\right)+M^8
 \end{align*}
 }
 {\tiny
 \begin{align*}
& A_4(L,M) =L^8 M^{16}+L^7 \left(-2 M^{24}+3 M^{22}-2 M^{20}+M^{18}+8 M^{16}+M^{14}-2 M^{12}+3
   M^{10}-2 M^8\right) \qquad \qquad \qquad \qquad \qquad \\
& +L^6
   \left(M^{32}-3 M^{30}+3 M^{28}-M^{26}-4 M^{24}-9
   M^{22}+15 M^{20}+13 M^{18}-2 M^{16}+13 M^{14} \right. \\
& \left.   +15 M^{12}-9 M^{10}-4 M^8-M^6+3
   M^4-3 M^2+1\right) \\
&  +L^5 \left(-4 M^{32}+16 M^{30}-14 M^{28}-21 M^{26}+24
   M^{24}+16 M^{22}-56 M^{20} \right. \\
& \left.   +17 M^{18}+100 M^{16}+17 M^{14}-56 M^{12}+16
   M^{10}+24 M^8-21 M^6-14 M^4+16 M^2-4\right)\\
&  +L^4
   \left(6 M^{32}-26 M^{30}+22
  M^{28}+40 M^{26}-59 M^{24}-64 M^{22}+82 M^{20} \right. \\
& \left.   +50 M^{18}-32 M^{16}+50 M^{14}+82
   M^{12}-64 M^{10}-59 M^8+40 M^6+22 M^4-26 M^2+6\right) \\
&  +L^3 \left(-4 M^{32}+16
   M^{30}-14 M^{28}-21 M^{26}+24 M^{24}+16 M^{22}-56 M^{20}+17 M^{18}+100
   M^{16}\right. \\
& \left.   +17 M^{14}-56 M^{12}+16 M^{10}+24 M^8-21 M^6-14 M^4+16 M^2-4\right)\\
& +L^2
   \left(M^{32}-3 M^{30}+3 M^{28}-M^{26}-4 M^{24}-9 M^{22}+15 M^{20}+13 M^{18}-2
   M^{16}+13 M^{14}+15 M^{12}-9 M^{10}\right. \\
& \left.   -4 M^8-M^6+3 M^4-3 M^2+1\right) \\
& +L \left(-2
   M^{24}+3 M^{22}-2 M^{20}+M^{18}+8 M^{16}+M^{14}-2 M^{12}+3 M^{10}-2
   M^8\right)+M^{16} 
 \end{align*}
 }  

\section{Appendix C} \label{app:C}

 {\tiny
 \begin{align*}
& p_2(\pm u) = \left(2 \left(L M^2+1\right)\right)^{-2}\left[\mp L
   \left(M^4-1\right)^2 \left(L+M^2\right)u \right.\\
&\left. +2 L^4 M^8+L^3 \left(-2 M^{12}+3 M^{10}+3 M^8+4 M^6+M^2-1\right)+L^2 \left(M^{12}-6
   M^{10}+5 M^8+12 M^6+5 M^4-6 M^2+1\right) \right.  \qquad \qquad \qquad \qquad \qquad\\
&   \left. +L \left(-M^{12}+M^{10}+4 M^6+3 M^4+3
   M^2-2\right)+2 M^4\right]
 \end{align*}  
 }
{\tiny
 \begin{align*}
& p_{-2}(\pm u) = (2 \left(L M^2+1\right))^{-\frac{3}{2}} \left[\pm L \left(M^4-1\right)^2 u \right.\\
&\left. +2 L^3 M^{10}-L^2 \left(M^{10}-5 M^8-2 M^6+M^2-1\right)+L \left(M^{10}-M^8+2
   M^4+5 M^2-1\right)+2\right]  \qquad \qquad \qquad \qquad \qquad  \qquad \qquad \qquad \qquad \qquad
\end{align*}
}
{\tiny
 \begin{align*}   
& p_{-4}(\pm u) = \left(2 \left(L M^2+1\right)\right)^{-7/2}
 \left[\pm L \left(M^4-1\right)^2 \left(L^4 \left(2 M^4-M^2+1\right) M^8+L^3 \left(-3
   M^{12}+9 M^{10}+M^8+M^6-M^4+2 M^2-1\right)  \right.\right. \\
&\left. \left. +L^2 \left(2 M^{12}-8 M^{10}+7 M^8+10
   M^6+7 M^4-8 M^2+2\right)+L \left(-M^{12}+2 M^{10}-M^8+M^6+M^4+9
  M^2-3\right)+M^4-M^2+2\right) u \right. \\
& \left.+2 L^7 M^{22}+L^6 \left(-4 M^{14}+15 M^{12}+6 M^{10}-5 M^8+3 M^4-2
   M^2+1\right) M^8 \right. \\
& \left.+L^5 \left(7 M^{22}-25 M^{20}+22 M^{18}+51 M^{16}-M^{14}-32
   M^{12}+13 M^{10}+10 M^8-2 M^6-3 M^4+3 M^2-1\right) \right. \\
& \left.+L^4 \left(-7 M^{22}+28
   M^{20}-21 M^{18}-31 M^{16}+62 M^{14}+82 M^{12}-50 M^{10}-30 M^8+33 M^6+18
   M^4-17 M^2+3\right) \right. \\
& \left.+L^3 \left(3 M^{22}-17 M^{20}+18 M^{18}+33 M^{16}-30
   M^{14}-50 M^{12}+82 M^{10}+62 M^8-31 M^6-21 M^4+28 M^2-7\right) \right.  \\
& \left.-L^2
   \left(M^{22}-3 M^{20}+3 M^{18}+2 M^{16}-10 M^{14}-13 M^{12}+32 M^{10}+M^8-51
   M^6-22 M^4+25 M^2-7\right) \right.  \\  
& \left. +L \left(M^{14}-2 M^{12}+3 M^{10}-5 M^6+6 M^4+15
   M^2-4\right)+2 \right]
 \end{align*}
 }
 
\section*{Acknowledgement}
This work was supported by Basic Science Research Program through the National Research Foundation of Korea (NRF) funded by the Ministry of Education, Science and Technology (No. NRF-2018R1A2B6005847). The second author was supported by 2018 Hongik University Research Fund.


\end{document}